\documentclass[reqno,draft]{amsart}

\usepackage{amsthm}
\usepackage{amssymb}
\usepackage{amsmath}

\newcommand{\labell}[1] {\label{#1}}

\newtheorem {Theorem}   {Theorem} 
\numberwithin{Theorem}{section}

\newtheorem {Lemma}[Theorem]    {Lemma}         
\newtheorem {Proposition}[Theorem]{Proposition}  
\theoremstyle{definition}

\theoremstyle{remark}
\newtheorem{Remark}[Theorem]{Remark}

\expandafter\chardef\csname pre amssym.def at\endcsname=\the\catcode`\@ 
\catcode`\@=11 
\def\undefine#1{\let#1\undefined} 
\def\newsymbol#1#2#3#4#5{\let\next@\relax 
 \ifnum#2=\@ne\let\next@\msafam@\else 
 \ifnum#2=\tw@\let\next@\msbfam@\fi\fi 
 \mathchardef#1="#3\next@#4#5}
\def\mathhexbox@#1#2#3{\relax 
 \ifmmode\mathpalette{}{\m@th\mathchar"#1#2#3}%
 \else\leavevmode\hbox{$\m@th\mathchar"#1#2#3$}\fi} 
\def\hexnumber@#1{\ifcase#1 0\or 1\or 2\or 3\or 4\or 5\or 6\or 7\or 8\or 
 9\or A\or B\or C\or D\or E\or F\fi} 
 
\font\teneufm=eufm10 
\font\seveneufm=eufm7
\font\fiveeufm=eufm5 
\newfam\eufmfam
\textfont\eufmfam=\teneufm 
\scriptfont\eufmfam=\seveneufm 
\scriptscriptfont\eufmfam=\fiveeufm 
 
\catcode`\@=\csname pre amssym.def at\endcsname 

\def    \eps    {\epsilon}
\def    \p      {\partial}
\def    \px     {\partial_x}
\def    \py     {\partial_y}
\def    \pt     {\partial_t}

\newcommand{\D}{{\mathfrak D}}

\def    \reals  {{\mathbb R}}
\def    \R      {{\mathbb R}}
\def    \integers       {{\mathbb Z}}

\def    \p      {\partial}

\def    \ssminus        {\smallsetminus}

\begin{document}


\setlength{\smallskipamount}{6pt}
\setlength{\medskipamount}{10pt}
\setlength{\bigskipamount}{16pt}





\title[A $C^2$-counterexample to the Hamiltonian Seifert Conjecture]{On 
the construction of a $C^2$-counterexample to the Hamiltonian Seifert 
Conjecture in $\R^4$}

\author[Viktor Ginzburg]{Viktor L. Ginzburg}
\author[Ba\c sak G\"urel]{Ba\c sak Z. G\"urel}

\address{Department of Mathematics, UC Santa Cruz, 
Santa Cruz, CA 95064, USA}
\email{ginzburg@math.ucsc.edu, basak@math.ucsc.edu}

\date{\today}

\thanks{The work is partially supported by the NSF and by the faculty
research funds of the University of California, Santa Cruz.}

\bigskip

\begin{abstract}
We outline the construction of a proper $C^2$-smooth function on $\R^4$ 
such that its Hamiltonian flow has no periodic orbits on at least one regular 
level set. This result can be viewed as a $C^2$-smooth counterexample to the 
Hamiltonian Seifert conjecture in dimension four. 
\end{abstract}

\maketitle

\section{Introduction} 
The goal of this notice is to outline the construction of a $C^2$-smooth 
proper function on $\R^4$ whose Hamiltonian flow has no periodic
orbits on at least one regular level set. A detailed construction of such a
function with complete proofs can be found in \cite{GG}. This result fits 
into the circle of questions generally referred to as the Seifert conjecture.

The original Seifert conjecture is the question whether or not there
exists a vector field on $S^3$ without zeros and periodic orbits.
In a broader context, the same question is asked for other manifolds or more
restricted classes of vector fields. Traditionally, examples of such vector
fields are called counterexamples to the Seifert conjecture. For closed 
manifolds of dimension greater than three with zero Euler characteristic,
$C^\infty$-smooth (in fact, real analytic) counterexamples were constructed 
by Wilson, \cite{wilson}. A $C^1$-smooth counterexample to the Seifert 
conjecture on $S^3$ was found by Schweitzer, \cite{schweitzer}, and 
a $C^1$-smooth volume--preserving counterexample on $S^3$ was constructed by 
G. Kuperberg, \cite{kug}. The ideas from both of these constructions play an 
important role in this paper. A real analytic counterexample on $S^3$ is due 
to K. Kuperberg, \cite{kugk,kuk}. The reader interested in a detailed survey 
of the results concerning the original Seifert conjecture 
should consult \cite{kuk:icm,kuk:notices}.

In a similar vein, the ``Hamiltonian Seifert conjecture'' is the question  
whether or not there exists a proper function on $\R^{2n}$ whose Hamiltonian 
flow has no periodic orbits on at least one regular level set. 
In dimensions greater than six, $C^\infty$-smooth counterexamples to 
the Hamiltonian Seifert conjecture were constructed 
by M. Herman, \cite{herman-fax,herman}, and simultaneously by one of 
the authors, \cite{gi:seifert}. In dimension six, a $C^{2+\eps}$-smooth 
counterexample was found by M. Herman, \cite{herman-fax,herman}. This
smoothness constraint was later relaxed to $C^\infty$  in 
\cite{gi:seifert97}. A very simple and elegant construction of new
$C^\infty$-smooth counterexamples in dimensions greater than four was
recently discovered by Kerman, \cite{ely:example}. We refer the reader
to \cite{gi:bayarea,gi:barcelona} for a detailed discussion of the 
Hamiltonian Seifert conjecture. 

An essential difference of the Hamiltonian case from the general one is 
manifested by the almost existence theorem, \cite{ho-ze:per-sol,ho-ze:book,str}, 
which asserts that almost all regular levels of a proper Hamiltonian have 
periodic orbits. To be more specific, for a $C^2$-smooth (and probably even 
$C^1$-smooth) proper Hamiltonian on $\R^{2n}$ regular values without periodic 
orbits form a zero measure set. Thus such regular levels are exceptional in the 
sense of measure theory.

\subsection*{Acknowledgments.} The authors are deeply grateful to 
Helmut Hofer, Anatole Katok, Ely Kerman, Krystyna Kuperberg, 
Mark Levi, Debra Lewis, Rafael de la Llave, Eric Matsui, and Maria Schonbek 
for useful discussions and suggestions.

\section{Main Results}
\labell{sec:main}

Recall that the characteristic foliation on a hypersurface $M$ in a
symplectic manifold $(W,\eta)$ is the one-dimensional foliation whose
leaves (called characteristics) are tangent to the field of directions
$\ker (\eta|_M)$.

Let $\R^{2n}$ be equipped with its standard symplectic structure.

\begin{Theorem}
\labell{thm:main}
There exists a $C^2$-smooth embedding $S^3\hookrightarrow\R^4$ such that 
its characteristic foliation has no closed characteristics. This embedding
can be chosen $C^0$-close and $C^2$-isotopic to an ellipsoid.
\end{Theorem}

As an immediate consequence we obtain

\begin{Theorem}
\labell{thm:main2}
There exists a proper $C^2$-function $F\colon\R^4\to \R$ such that
the level $\{F=1\}$ is regular and has no periodic orbits. In addition, 
$F$ can be chosen so that this level is $C^0$-close and $C^2$-isotopic to an
ellipsoid. 
\end{Theorem}

Note that since all values of $F$ near $F=1$ are regular, by the 
almost existence theorem, almost all levels of
$F$ near this level carry periodic orbits. 

\begin{Remark}
It is likely that our construction of the embedding 
$S^3\hookrightarrow\R^4$ gives in fact a $C^{2+\alpha}$-embedding without
closed characteristics. However, checking this requires considerably more
work than needed to verify $C^2$-smoothness.
\end{Remark}

\begin{Remark}
\labell{rmk:main}
Similarly to its higher-dimensional counterparts, \cite{gi:seifert,gi:seifert97},
Theorem \ref{thm:main} extends to other symplectic manifolds as follows.
Let $(W,\eta)$ be a four-dimensional symplectic manifold and let
$i\colon M\hookrightarrow W$ be a $C^\infty$-smooth embedding such that 
$i^*\eta$ has only a 
finite number of closed characteristics. Then there exists a $C^2$-smooth 
embedding $i'\colon M\hookrightarrow W$, which is $C^0$-close and isotopic 
to $i$, such that ${i'}^*\eta$ has no closed characteristics.
\end{Remark}

In what follows we outline the proof of Theorem \ref{thm:main}. 
The idea of the proof is to adjust Schweitzer's construction, \cite{schweitzer}, 
of a $C^1$-flow on $S^3$ without periodic orbits to make it embeddable into 
$\R^4$ as a Hamiltonian flow. This is done by constructing a version of 
Schweitzer's plug with suitable properties similar to those of G. Kuperberg's 
plug, \cite{kug}. The essential feature of the construction is that Schweitzer's 
plug has to be built using a ``very smooth'' Denjoy flow 
(see Remark \ref {rmk:smooth}).

\section{Construction of symplectic embedding}
\labell{sec:embedding}

Let us first fix some notations. Throughout this paper $\sigma$ 
denotes the standard symplectic form on $\R^{2m}$ or the pull-back of
this form to $\R^{2m+1}$ by the projection $\R^{2m+1}\to \R^{2m}$ along
the first coordinate; $I^{2m}$ stands for a cube in $\R^{2m}$ whose
edges are parallel to the coordinate axis. The product $[a,b]\times I^{2m}$ is
always assumed to be embedded into $\R^{2m+1}$ (henceforth, the standard 
embedding) so that the interval $[a,b]$ is 
parallel to the first coordinate. We refer to the direction along the
first coordinate $t$ (time) in $\R^{2m+1}$ (or $[a,b]$ in 
$[a,b]\times I^{2m}$) as the vertical direction.
All maps whose smoothness is not specified are $C^\infty$-smooth.

Theorem \ref{thm:main}, as do similar theorems in dimensions greater
than four, follows from the existence of a symplectic plug. The definitions
of a plug vary considerably (see \cite{gi:seifert,ely:example,kug}), and here 
we use the one more suitable for our purposes.

A \emph{$C^k$-smooth symplectic plug} in dimension $2n$ is a $C^k$-embedding $J$ of 
$P=[a,b]\times I^{2n-2}$ into 
$P\times \R\subset\R^{2n}$ such that

\begin{enumerate}

\item[P1.] 
\emph{The boundary condition}: The embedding $J$ is the identity
embedding of $P$ into $\R^{2n-1}$ near the boundary $\partial P$. Thus
the characteristics of $J^*\sigma$ are parallel to the vertical direction
near $\partial P$.

\item[P2.]
\emph{Aperiodicity}:
The characteristic foliation of $J^*\sigma$ is aperiodic, i.e., has no 
closed characteristics.

\item[P3.] 
\emph{Existence of trapped trajectories}:
There is a characteristic of $J^*\sigma$ beginning on $\{a\}\times I^{2n-2}$ 
that never exits the plug. Such a characteristic is said to be trapped in $P$.

\item[P4.]
The embedding $J$ is $C^0$-close to the standard embedding and $C^k$-isotopic
to it.

\item[P5.]
\emph{Matched ends condition}:
For a characteristic that meets both the bottom and the top of the plug,
its top end lies exactly above the bottom point .

\end{enumerate}

\begin{Theorem}
\labell{thm:plug}
In dimension four, there exists a $C^2$-smooth symplectic plug.
\end{Theorem}

Theorem \ref{thm:main} readily follows from Theorem \ref{thm:plug} by
applying a symplectic plug to perturb an irrational ellipsoid in $\R^4$.

\begin{proof}[Proof of Theorem \ref{thm:plug}]
By the standard symmetry argument it suffices to construct a semi-plug, 
i.e., ``plug'' satisfying only conditions (P1)-(P4). 
We will construct a semi-plug by perturbing the natural embedding 
of $[a,b]\times I^2$ into $\R^4$.

Let $\Sigma$ be the torus $T^2$ with coordinates $(x,y)$, punctured at 
$(x_0,y_0)$, i.e., with
a neighborhood of $(x_0,y_0)$ deleted. There exists a symplectic bridge 
immersion of $(\Sigma, dx\wedge dy)$ into some cube $I^2$ with the standard 
symplectic structure. Hence, there exists an embedding, referred to 
in what follows as the standard embedding,
$$
M=[-1,1]\times \Sigma \hookrightarrow [a,b]\times I^2\subset \R^3
\subset \R^4
$$
such that the pull back of $\sigma$ is $dx\wedge dy$. From now on, we will
identify $M$ with the image of this embedding. 

The embedding $J$ is obtained from the standard embedding by altering it within
$M$. This is done in a few steps. First, we consider an embedding of $M$ into 
some four-dimensional symplectic manifold $(W, \sigma_W)$ such that
the pull-back of $\sigma_W$ is still $dx\wedge dy$. 
Then we $C^0$-perturb this embedding so that the characteristic vector field of
the new pull-back will have the properties similar to those of Schweitzer's
plug. By the symplectic neighborhood theorem, a neighborhood of $M$ in
$W$ is symplectomorphic to that of $M$ in $\R^4$. This will allow us to
turn the embedding $M\hookrightarrow W$ into the required embedding
$J\colon M\hookrightarrow \R^4$.

As in Schweitzer's construction, the characteristic vector field of $J$ 
is a modification of a Denjoy vector field on $T^2$. This modification 
is obtained by extending the Denjoy vector field to $[-1,1]\times T^2$
and then restricting it to $\varphi(M)$ for a suitable embedding 
$\varphi\colon M\to [-1,1]\times T^2$. Likewise, $J$ is essentially equal to
$j\varphi$ for an embedding $j$ of $[-1,1]\times T^2$ into 
$W$. Hence we start with defining $j$ on 
$[-1, 1]\times T^2$.

Let $W=(-2,2)\times S^1\times T^2$ with coordinates
$(t, x,u,y)$ and symplectic form 
$\sigma_W=dt\wedge dx+du\wedge dy$ and let
$$
j_0\colon [-1,1]\times T^2\hookrightarrow W;
\quad j_0(t,x,y)=(t,x,x,y).
$$
Consider a mapping $K\colon [-1,1]\times T^2\to S^1$ to be specified later on
and define
$$
j\colon [-1,1]\times T^2\to (-2,2)\times S^1\times T^2
\quad\text{by}\quad j(t,x,y)=(t,x,K,y).
$$
(To illustrate this definition note that when $K(t,x,y)=x$, the embedding $j$ 
turns into $j_0$.)
The characteristic vector field of the pull-back form $j^*\sigma_W$ is 
$$
v=(\px K)\pt-(\pt K) \px + \py.
$$
\begin{Remark}
\labell{rmk:j}
To justify the definition of $j$, let us view
the annulus $[-1,1]\times S^1$ as the symplectic manifold with symplectic
form $dt\wedge dx$ and the product $[-1,1]\times T^2$ as the extended
phase space with the $y$-coordinate being the time-variable. Then
we can regard $K$ as a (multi-valued) time-dependent Hamiltonian on 
$[-1,1]\times S^1$ and view $W$ as the further extended time-energy
phase space with the cyclic energy-coordinate $u$. Then $j$ is the graph
of the time-dependent Hamiltonian $K$ in the extended time-energy phase
space $W$. Now it is clear that $v$ is just the Hamiltonian vector field 
of $K$.
\end{Remark}

In what follows we choose $\varphi\colon M\to [-1,1]\times T^2$
so as to make a neighborhood $U$ of $j_0\varphi(M)$ in $W$ symplectomorphic
to a neighborhood of $M$ in $\R^4$. Then, if $j$ takes values in $U$ and 
matches $j_0$ near the boundary $\p M$, we can turn $j\varphi$ into
the desired embedding $J$.  However, to ensure that
(P1)-(P4) are satisfied we need to impose some requirements on $K$.
Let us now describe these requirements.

Let $\py+h\px$ be a Denjoy vector field on $T^2$ to be specified 
in Section \ref{sec:pfK1}
and let $\D$ be the Denjoy continuum for this field. (See, e.g., 
\cite{KH,schweitzer} for an introduction to Denjoy maps and vector fields.) 
Here we only point out that $(x_0,y_0)$ is chosen in the complement of $\D$.
Fix a small neighborhood $V$ of $(x_0,y_0)$ disjoint from $\D$. (This
neighborhood will contain the puncture, i.e., $\p \Sigma \subset V$.)
Consider the tubular neighborhood of the line $(t,x_0,y_0+t)$
in $[-1,1]\times T^2$ of the form 
$\{ (t,x,y+t)\mid (x,y)\in V,~ t\in [-1,1]\}$.
Fix also a small neighborhood of the boundary
$\p ([-1,1]\times T^2)$ and denote by $N$ the union of these neighborhoods.

\begin{Proposition}
\labell{prop:K1}
There exists a $C^2$-smooth mapping $K\colon [-1,1]\times T^2\to S^1$ such
that 

\begin{enumerate}

\item[K1.] $v$ is equal to the Denjoy vector field (i.e.,
$\px K =0$ and $\pt K=-h$) at every point of $\{0\}\times \D$;

\item[K2.] the $t$-component of $v$ is positive (i.e., $\px K> 0$) 
on the complement of $\{0\}\times \D$;

\item[K3.]$K$ is $C^0$-close\footnote{More specifically, 
for any $\eps>0$ there exists a $K$ satisfying (K1)-(K2) and (K4) such that
$\parallel K-K_0\parallel <\eps$. The required value of $\eps$ is 
determined by the size of the neighborhood $U$ in the symplectic 
neighborhood theorem; see below.}  to the map 
$K_0\colon (t,x,y)\mapsto x$;

\item[K4.]$K=K_0$ on $N$.
\end{enumerate}
\end{Proposition}

Let us defer the proof of this proposition to Section \ref{sec:pfK1} and 
finish the proof of Theorem \ref{thm:plug}. From now on we assume that $K$ 
is as in Proposition \ref{prop:K1}.

By (K1) and (K2), $v$ has a trapped trajectory and is aperiodic.
Indeed, by (K1), $\{0\}\times \D$ is invariant under the flow of $v$ and
on this set the flow is a Denjoy flow. By (K2), the vertical 
component of $v$ is non-zero unless the point is in $\{0\}\times \D$. 
This implies that periodic orbits can only occur within $\{0\}\times \D$.
Since the Denjoy flow is aperiodic, so is the entire flow of $v$.
Furthermore, it is easy to see that since $\{0\}\times \D$ is invariant, there
must be a trapped trajectory. 

Let
$$
\varphi\colon M=[-1,1]\times \Sigma \to [-1,1]\times T^2;
\quad \varphi(t,x,y)=(t,x,y+t).
$$
Then $(j_0\varphi)^*\sigma_W= dx\wedge dy$.
The argument similar to the proof of 
the symplectic neighborhood theorem, \cite[Lemma 3.14]{mcduff-sal},
shows that a ``neighborhood'' of $M$ in $\R^4$ is
symplectomorphic to a ``neighborhood'' $U$ of $j_0\varphi(M)$ 
in $W$ (see \cite[Section 4]{gi:seifert} and \cite{GG} for details). More 
precisely, for a small $\delta>0$, there exists a symplectomorphism
$$
\psi\colon M\times (-\delta,\delta)\to U\subset W
$$
extending $j_0\varphi$, i.e., such that $\psi|_{M}=j_0\varphi$.

By (K3),  $j$ is $C^0$-close to $j_0$ 
and $j=j_0$ on $N$ by (K4). Hence, $j$ can be assumed to take values
in $U$. Set
$$
J=\psi^{-1}j\varphi
$$
on $M$ and extend $J$ as the standard embedding to $[a,b]\times I^2\ssminus M$.
Then $(J^*\sigma)|_M=(j\varphi)^*\sigma_W$.

The characteristic vector field of $J^*\sigma$ is $\pt$ in the
complement of $M$ and $(\varphi^{-1})_*v$ on $M$. Since
$(\varphi^{-1})_*v=\pt$ near $\p M$, these vector fields match smoothly
at $\p M$. Clearly, (P1) is satisfied. Since $v$ has a trapped trajectory and 
is aperiodic, the same is true for $(\varphi^{-1})_*v$, i.e., 
the conditions (P2) and (P3) are met.  The condition (P4) is easy to 
verify. Hence, $J$ is a semi-plug.

\end{proof}

\begin{Remark}
\labell{rmk:measure-zero}
In the proof of Proposition \ref{prop:K1} we will not require the
Denjoy continuum $\D$ to have zero measure. As a consequence, the
union of characteristics entirely contained in the semi-plug 
can have Hausdorff dimension two because this set is the 
image of $\D$ by a $C^2$-smooth embedding.
\end{Remark}

\section{Construction of the function $K$}
\labell{sec:pfK1}

\subsection{Denjoy vector field}
The first step in the construction of $K$ is to choose
an appropriate Denjoy vector field on $T^2$. 

\begin{Lemma}
\labell{lemma:main}
There exists a Denjoy vector field $\py +h\px$ which 
is $C^{1+\alpha}$ for all $\alpha \in (0,1)$ and such that
\begin{enumerate}
\item[D1.] $\px h$ vanishes on the Denjoy continuum $\D$;

\item[D2.] $\int_0^x (\px h(\xi,y))^2\,d\xi$ is $C^2$ in $(x,y)$.
\end{enumerate}
\end{Lemma}

\begin{proof}[Sketch of the proof of Lemma \ref{lemma:main}]
We start with the Denjoy map $\Phi\colon S^1\to S^1$ defined using 
intervals $I_n$ of length
\begin{equation}
\labell{defn:$l_n$}
l_n:=k_\beta(|n|+2)^{-1}(\log\,(|n|+2))^{-1/\beta}
\end{equation}
for some $\beta\in (0,1)$.  Here $k_\beta$ is 
a constant depending on $\beta$ chosen so that 
$\sum_{n\in\integers}\,l_n< 1$. (See \cite{KH} for details.) 
The importance of our choice of $l_n$ is 
that the series $\sum_{n\in\integers}\,l_n$ converges very slowly, which
results in a small Denjoy continuum $\D_0=S^1 \ssminus 
\bigcup_{n\in \integers}\,Int(I_n)$.  
This slow convergence is the main factor which ensures 
that $\Phi$ is sufficiently smooth. More specifically, it is well
known and easy to see that  $\Phi$ is
$C^{1+\alpha}$ for any $\alpha\in (0,1)$ regardless of what $\beta$
we have taken. Furthermore, we claim that
\begin{equation}
\labell{eq:Phi}
(\Phi'-1)^2\quad\text{is}\quad C^1.
\end{equation}
Indeed, $\Phi'-1$ vanishes on $\D_0$. Then the fact that $\Phi'-1$ is
$C^\alpha$ with $\alpha>1/2$ readily implies that $(\Phi'-1)^2$ is 
everywhere differentiable and that its derivative vanishes along $\D_0$.
The proof that this derivative is continuous relies explicitly
on \eqref{defn:$l_n$}. Using \eqref{defn:$l_n$}, one can show
that
$$
\Bigl\| \frac{d}{dx}(\Phi'-1)^2\Big|_{I_n}\Bigr\|_{\infty}
\to 0\quad\text{as}\quad n\to\infty,
$$
which ensures the continuity of the derivative.

Let now $\px+h\py$ be a suitably defined Denjoy vector field 
on $T^2$ for $\Phi$. The $C^{1+\alpha}$-smoothness of $h$ and (D1) do not 
present a problem. However, the verification of (D2) is more involved.
On the conceptual level, the proof is similar to that of \eqref{eq:Phi}
outlined above. To illustrate this point, observe that $\Phi'-1$
can be taken as an analogue of $\px h$ and one may expect these
two functions to have the same smoothness properties. (In fact, $(\px h)^2$
is $C^1$ in the $x$-variable.) The actual proof of (D2) requires a direct
verification of the existence and continuity of derivatives using their
asymptotic behavior on the intervals $I_n$; see \cite{GG} for details.
\end{proof}

\begin{Remark}
\labell{rmk:smooth}
The Denjoy map defined by \eqref{defn:$l_n$} is essentially as smooth
as a Denjoy map can be, up to using functions growing slower than 
logarithms, e.g., iterations of logarithms. The next significant 
improvement in smoothness would be to have $\log \Phi'$ of bounded variation
or satisfying the Zygmund condition which is impossible; see \cite{HS,KH}.
\end{Remark}

\subsection{Construction of K}
First we extend $h$ to a $C^1$-function on 
$H\colon [-1,1]\times T^2\to \R$ so that 
\begin{equation}
\labell{eq:H2}
\px H (t,x,y)= \px h(x,y)+o(t)\quad\text{uniformly in $(x,y)$}
\end{equation}
and the function
\begin{equation}
\labell{eq:H3}
\int_0^t H(\tau,x,y)\,d\tau\quad\text{is $C^2$ in $(t,x,y)$}.
\end{equation}
The value $H(t,x,y)$ is obtained by averaging $h$ over the 
square with side $t^s$, for some $s>0$, centered at $(x,y)$.
(This averaging argument is somewhat similar to the one from \cite{kug}.) 
Then using only the fact that $h$ is $C^{1+\alpha}$ with $\alpha>1/2$, it is 
not hard to show that $H$ satisfies \eqref{eq:H2} and \eqref{eq:H3}.

Now we are in a position to define the function $K$. We focus on 
a small neighborhood $U$ of $\D$ in $T^2$ and small values of $t$.
In other words, $K$ will be defined only on $[-\eps, \eps]\times U$,
where $\eps>0$ is sufficiently small. Furthermore, we consider only
the requirements (K1) and (K2). For this is where the essence of the
problem lies, whereas the conditions (K3) and (K2) are purely technical.
With these restrictions on the domain of $K$ and the hypotheses this
function should satisfy, we can take $K$ to be real-valued rather than 
circle-valued. (We refer the reader to \cite{GG} for a complete construction 
of $K$ and a detailed proof of (K1)-(K4).)

We will construct the function $K$ of the form
\begin{equation}
\labell{eq:K}
K(t,x,y)=-\int_0^t H(\tau,x,y)\,d\tau + A(x,y)+f(x,y)t^2.
\end{equation}
Here the ``constant'' of integration $A$ is a $C^2$ function $U\to \R$
such that
\begin{equation}
\labell{eq:A1}
\px A\geq 2(\px h)^2\quad\text{and}\quad \px A(x,y)=0 \quad\text{iff}\quad
(x,y)\in \D.
\end{equation}
The existence of such a function depends heavily on (D1) and (D2). Indeed,
$\px A$ is $C^1$, and $\px A$ and $(\px h)^2$ both vanish on $\D$. Thus,
\eqref{eq:A1} forces $(\px h)^2$ to be at least differentiable in $x$ at the
points of $\D$. Conversely, using (D1) and (D2), we can set
$$
A(x,y)= 2\int_0^x (\px h(\xi,y))^2\,d\xi+ a(x,y),
$$
where $a$ is a $C^\infty$-function on $U$ such that $\px a\geq 0$ 
and $\px a=0$ exactly on $\D$. Then $A$ will satisfy \eqref{eq:A1}.

The correction function $f\colon U\to \R$ is $C^\infty$-smooth 
and chosen so that
$$
\px f=2
$$
(One can show that such a function does exist.)

Let us now verify the properties of $K$. First note that $K$ is indeed
$C^2$-smooth by \eqref{eq:H3}. Condition (K1) is also obvious: 
since $H|_{t=0}=h$, we have $\pt K|_{t=0} = -H|_{t=0}= -h$ and
$\px K|_{\{0\}\times \D} = \px A|_{\D}= 0$ by \eqref{eq:A1}.

Let us now turn to (K2). We will first show that
\begin{equation}
\labell{eq:inequality}
\px K= -\int_0^t\px H\,d\tau +\px A +2t^2   \geq 0
\end{equation}
and then prove that the equality occurs only on $\{0\}\times \D$.

By \eqref{eq:H2}, we have
$$
\px K= \px A -t\px h  +t^2
+\left(t^2 + o(t^2)\right).
$$
Clearly,
\begin{equation}
\labell{eq:go4}
t^2 + o(t^2) \geq 0
\end{equation}
if $\eps>0$ is small. Hence, to verify \eqref{eq:inequality} it suffices 
to show that
\begin{equation}
\labell{eq:pfk3}
\px A -t\px h  +t^2 \geq 0.
\end{equation}
By \eqref{eq:A1}, this inequality is a consequence of
$$
2 (\px h)^2
-t\px h  +t^2  \geq 0.
$$
Here all the terms are non-negative except, maybe, $-t\px h$.
Hence, it suffices to show that at least one of the following two 
inequalities holds:
\begin{eqnarray}
\labell{eq:ineq1}
2(\px h)^2- t\px h  &\geq& 0,\\
\labell{eq:ineq2}
-t\px h +t^2 &\geq& 0.
\end{eqnarray}
Inequality \eqref{eq:ineq1} holds if (but not only if)
$$
|t|\leq 2 |\px h|
$$
and \eqref{eq:ineq2} holds if (but not only if)
$$
|t|\geq |\px h|
$$
It is clear that at least one of these inequalities 
holds. This proves \eqref{eq:inequality}.

To finish the proof of (K2) we need to show that the equality
in \eqref{eq:inequality} implies that $t=0$ and $(x,y)\in 
\D$. Thus, assume that $\px K(t,x,y)=0$. Then 
\eqref{eq:go4} and \eqref{eq:pfk3} must become equalities. 
The equality \eqref{eq:go4}
is possible only when $t=0$. Setting $t=0$ in 
the equality \eqref{eq:pfk3}, we conclude that $\px A (x,y)=0$ and hence 
$(x,y)\in \D$ by \eqref{eq:A1}. This concludes the proof of (K2).

\begin{Remark}
The actual construction of $K$, and, in particular, of $A$ and $f$, given
in \cite{GG} differs considerably from the one above. In fact, when $K$ is
defined on $[-\eps,\eps]\times T^2$, the constant of integration $A$ is
a function $T^2\to S^1$. Moreover, some extra care is needed in the definition
of $f$. However, these difficulties are rather technical and our description
captures well the main idea of the definition of $K$. 
\end{Remark}

\end{document}